**УДК 372.8+519.6**

**Словак Катерина Іванівна**, кандидат педагогічних наук, старший викладач кафедри вищої математики, Криворізький економічний інститут ДВНЗ «Криворізький національний університет», м. Кривий Ріг, e-mail: Slovak_kat@mail.ru

# МЕТОДИКА ПОБУДОВИ ОКРЕМИХ КОМПОНЕНТІВ МОБІЛЬНОГО МАТЕМАТИЧНОГО СЕРЕДОВИЩА «ВИЩА МАТЕМАТИКА»

**Анотація**

Актуальність матеріалу, викладеного у статті, обумовлена необхідністю розробки й упровадження у процес навчання високотехнологічного інформаційно-комунікаційного освітньо-наукового середовища. У роботі розглянуто один із прикладів такого web-орієнтованого середовища для навчання математичних дисциплін студентів ВНЗ – мобільне математичне середовище. На прикладі ММС «Вища математика» продемонстровано основні технології і засоби для побудови таких компонентів середовища, як лекційні демонстрації, динамічні моделі, тренажери, генератори навчальних завдань, навчальні експертні системи, індивідуальні домашні завдання, приклади розв'язування завдань, задачі для практичних занять тощо.

**Постановка проблеми**. На процес розвитку і формування майбутнього фахівця з вищою освітою найбільше впливає середовище, у якому він навчається. Саме тому, важливою проблемою для ВНЗ є розробка й упровадження у навчальний процес високотехнологічного інформаційно-комунікаційного освітньо-наукового середовища. Таке середовище має сприяти активізації навчальної діяльності студентів, розкриттю їхнього творчого потенціалу, збільшенню ролі самостійної роботи, розвитку професійної мобільності, формуванню фахових компетентностей, підвищенню рівня інформаційно-комунікаційної підготовки, а також відповідати потребам інформаційного суспільства і світовим освітнім стандартам. Прикладом такого середовища для навчання математичних дисциплін студентів ВНЗ є web-орієнтоване математичне середовище – *мобільне математичне середовище (ММС).*

**Аналіз останніх досліджень**. Як зазначають С. О. Семеріков, К. І. Словак

[1−2], під час побудови такого середовища (як ядра) доцільно використовувати вільно поширювані web-орієнтовані системи комп’ютерної математики (СКМ), що інтегрують в собі послуги різних систем за допомогою клієнт-серверних технологій і таких засобів ІКТ навчання математики, як мультимедійні демонстрації, динамічні математичні моделі, тренажери та експертні системи навчального призначення. Отже, основними складовими ММС є *обчислювальне ядро* й *інформаційне забезпечення*, що містить навчально-методичні і допоміжні інформаційні матеріали. Докладну характеристику складових, структури ММС та приклади використання розробленого ММС «Вища математика» [3], призначеного для підтримки навчання вищої математики студентів економічних спеціальностей, розглянуто в роботах [1–7].

**Метою статті** є висвітлення принципів створення та технології побудови окремих компонентів інформаційного забезпечення ММС «Вища математика».

**Основна частина**. Використання таких складових ММС, як лекційні демонстрації, динамічні моделі, тренажери та генератори навчальних завдань, передбачає багаторазове виконання обчислень для різних значень вхідних параметрів, тому під час їх розробки доцільно використати візуальні елементи управління типу «поле для введення», «повзунок», «прапорець», «меню вибору», для створення яких використовують відповідні функції обчислювальної складової ММС – Web-СКМ Sage.

Визначення кожного елементу управління здійснюється мовою Python за допомогою декоратора `@interact`, після якого ключовим словом `def` оголошується сама функція та її ім’я `name`. Якщо `name` не використовується в подальших розрахунках, то використовують символ «_» – нижнє підкреслення, що позначає безіменну функцію. Потім, деякій змінній `a` присвоюємо результат виконання функції, що відповідає потрібному елементу управління. Для ілюстрації зовнішнього вигляду створених елементів управління необхідно використати функцію `show()`, що показує як сам елемент управління, так і поточне значення змінної `a`, що виводиться у полі графічних побудов. Детальний опис різних функцій для програмування елементів управління подано у табл. 1.

Зазвичай, створення комп’ютерних моделей потребує одночасного виведення на екран декількох елементів управління, раціональне й естетичне розташування яких досягається за рахунок використання опції `layout` з відповідними ключами.

*Таблиця 1*

**Відомості щодо створення елементів управління**

**«Повзунок 1»**

***Функція***: `slider(vmin,vmax,step_size,default,label,display_value)`,

`vmin` – основний параметр для задання мінімального значення;

`vmax` – основний параметр для задання максимального значення;

`step_size` – додатковий параметр для задання кроку зміни числових значень;

`default` – додатковий параметр для задання значення за замовчуванням;

`label` – додатковий параметр для задання надпису ліворуч від елемента;

`display_value` – додатковий параметр логічного типу для регулювання виведенням на екран поточного значення.

***Приклади:***

```
@interact
def name(a = slider(1,9,1,default=4,label="α",)):
   show(a)
```

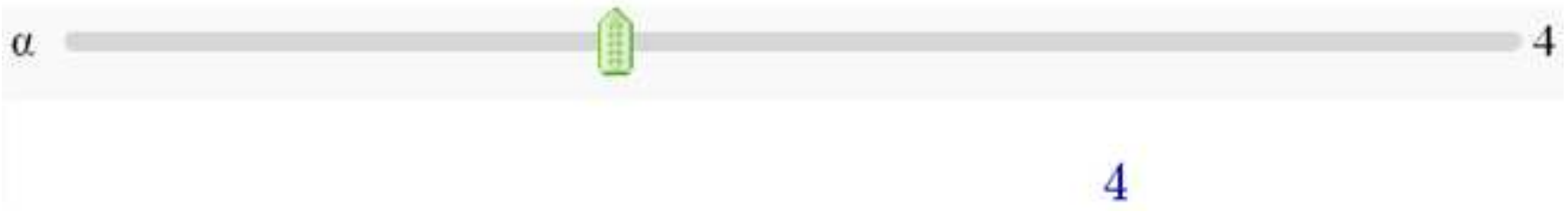

Значеннями «повзунка» можуть бути елементи списку різноманітної природи: `[1..100]` – цілі числа від 1 до 100; `[1, 'x', 'abc', 2/3]` – 4 елементи різної природи.

**«Повзунок 2»**

***Функція:*** `range_slider(vmin,vmax,step_size,default,label])`,

`vmin` – основний параметр для задання мінімального значення;

`vmax` – основний параметр для задання максимального значення;

`step_size` – додатковий параметр для задання кроку зміни числових значень;

`default` – додатковий параметр для задання значень за замовчуванням у форматі (`value_left,value_right`);

`label` – додатковий параметр для задання надпису ліворуч від елемента.

***Приклади:***

```
@interact
def _(a = range_slider (1,10,1, default=(4,5),
```

```
label = 'Діапазон' )):
    show(a)
```

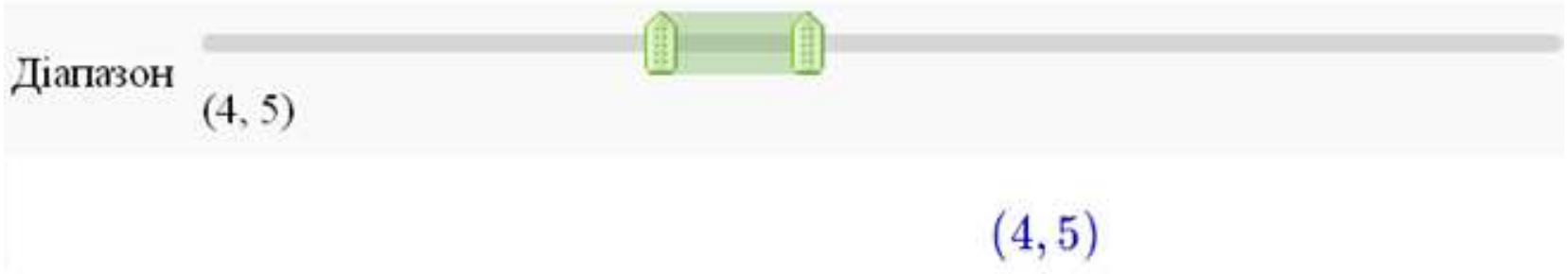

**«Прапорець»**

***Функція***: `checkbox(default,label)`,

`default` – основний параметр для задання стану «прапорця» за замовчуванням (набуває значення `false` або `true`);

`label` – додатковий параметр для задання надпису ліворуч від елемента.

***Приклади:***

```
@interact
def _(a=checkbox(False,"Показувати відповідь")):
    show(a)
```

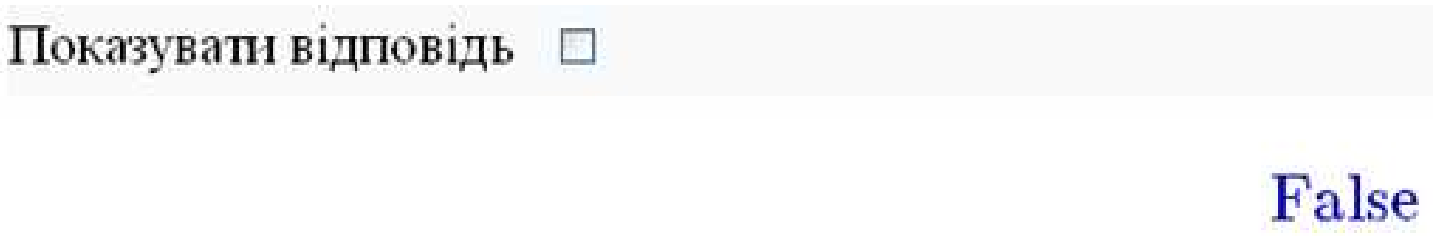

```
@interact
def _(a=checkbox(True,"Показувати відповідь")):
    show(a)
```

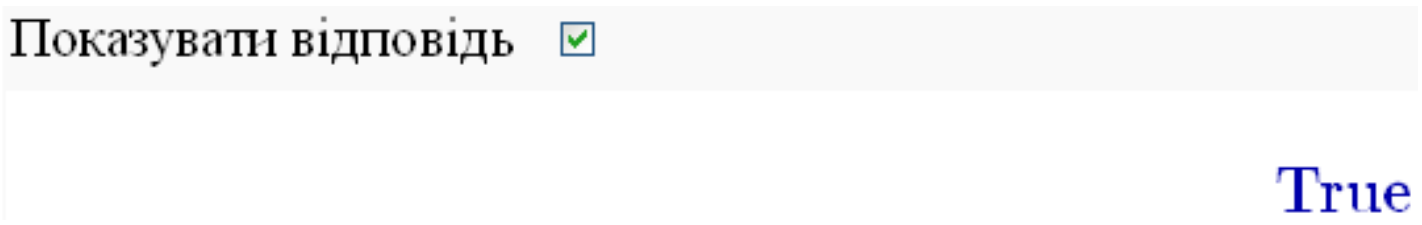

**«Меню вибору»**

***Функція***: `selector(values,label,default,nrows,ncols,width, buttons)`,

`values` – основний параметр для задання значення пунктів меню вибору, що можуть зазначатися переліком елементів – `[val1,val2,val3,...]` або діапазоном елементів – `[val_start..val_finish]`;

`label` – додатковий параметр для задання надпису ліворуч від елемента;

`default` – додатковий параметр для задання значень за замовчуванням;

`nrows` – додатковий параметр для задання кількості рядків у поданні пунктів меню вибору (при поданні пунктів меню вибору у вигляді кнопок);

`ncols` – додатковий параметр для задання кількості стовпчиків у поданні пунктів меню вибору (у разі подання пунктів меню вибору у вигляді кнопок);

`width` – додатковий параметр для задання ширини кнопок (у разі подання пунктів меню вибору у вигляді кнопок);

`buttons` – додатковий параметр логічного типу: при встановленому значенні `true` меню вибору подається у вигляді кнопок, при встановленому значенні `false` (за замовчуванням) – у вигляді списку, що розкривається.

***Приклади:***

```
@interact
def  _(a=selector([1..5],"Виберіть  значення",  default=2,
buttons=false)):
    show(a)
```

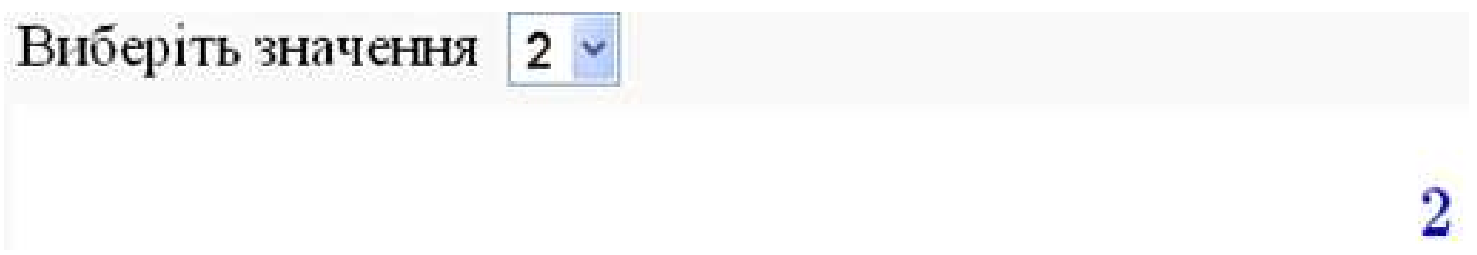

```
@interact
def  _(a=selector(['Приклад  1','Приклад  2',  'Приклад  3',
'Приклад   4',   'Приклад   5',   'Приклад   6'],label="",
default='Приклад 6', nrows=3, ncols=2, width=15 )):
    show(a)
```

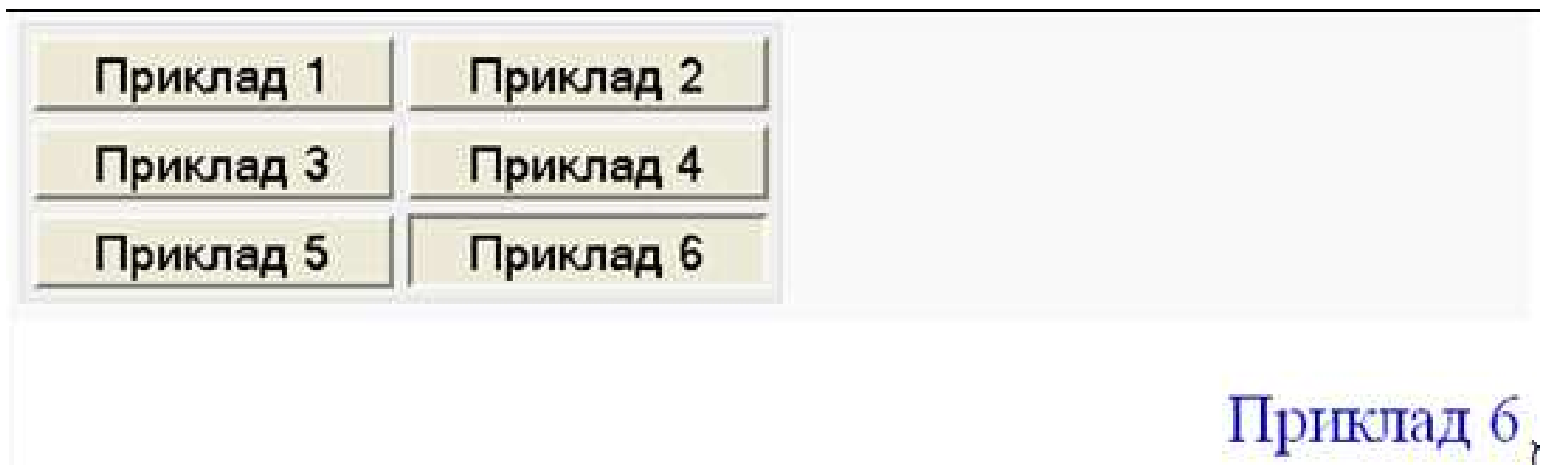

**«Поле для введення»**

***Функція:*** `input_box(default,label,type,width,kwargs)`,

`default` – основний параметр для задання значення, що повертається функцією за замовчуванням;

`label` – додатковий параметр для задання надпису ліворуч від елемента;

`type` – додатковий параметр для визначення типу даних, що вводяться;

`width` – додатковий параметр для задання ширини поля;

`kwargs` – додатковий параметр для підключення одного з існуючих словників.

***Приклади:***

```
@interact
def _(a=input_box("2+89", 'Введіть значення', width=10)):
    show(a)
```

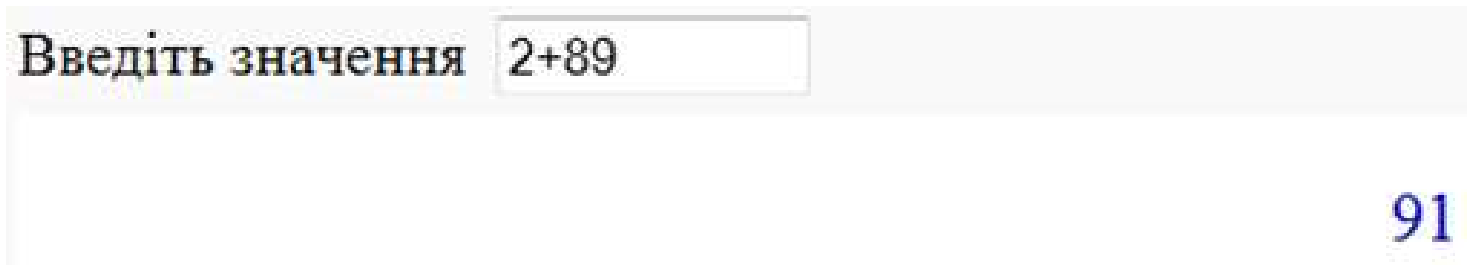

```
@interact
def _(a=input_box('Sage', label="", type=str)):
    show(a)
```

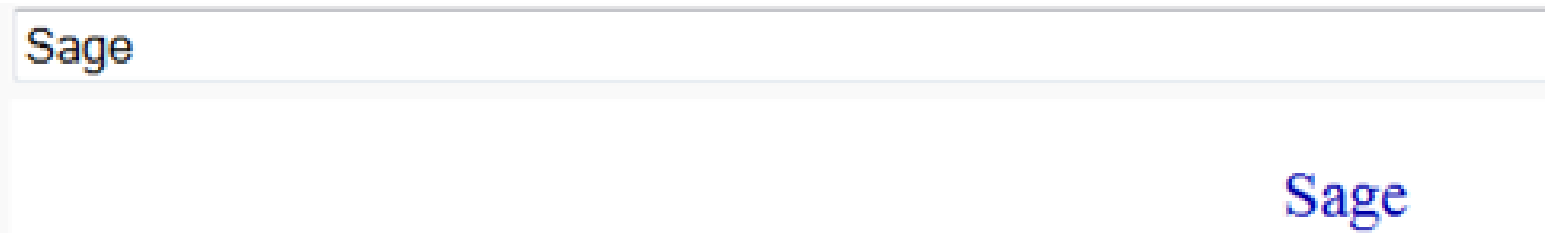

**«Комірки для введення»**

***Функція:***

```
input_grid(nrows,ncols,default,label,to_value=lambda
x:x,width),
```

`nrows` – основний параметр для задання кількості рядків;

`ncols` – основний параметр для задання кількості стовпчиків;

`default` – основний параметр для задання початкових значень у комірках;

`label` – додатковий параметр для задання надпису ліворуч від елемента;

`to_value=lambda x:x` – основний параметр для формування і виведення заданих даних у вигляді списку;

`width` – додатковий параметр для задання ширини комірок.

***Приклади:***

```
@interact
def _(a=input_grid(3,4,default=[1,2,3,4,5,6,7,8,9,
10,11,12], label='Матриця A', to_value=lambda x:x,
width=2)):
    show(a)
```

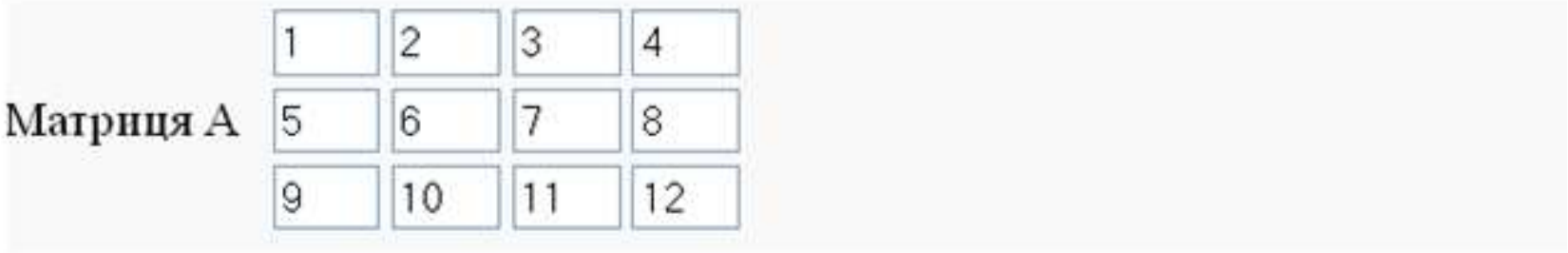

$[[1,2,3,4],[5,6,7,8],[9,10,11,12]]$

**«Поле вибору кольору»**

***Функція:*** `color_selector(default,label,widget,hide_box)`,

`default` – основний параметр для задання кольору у палітрі `(R,G,B)`;

`label` – додатковий параметр для задання надпису ліворуч від елемента;

`widget` – основний параметр для задання вигляду діалогового вікна, за замовчуванням присвоюється значення `jpicker`, також може набувати значень `farbtastic` або `colorpicker`.

`hide_box` – основний параметр для відображення вікна введення кольору у шістнадцятковому форматі, за замовчуванням присвоюється значення False.

***Приклади:***

```
@interact
def _(c = color_selector((1,0,0), "Колір", widget=
'farbtastic', hide_box=False)):
    show(plot(sin(7*x), color = c))
```

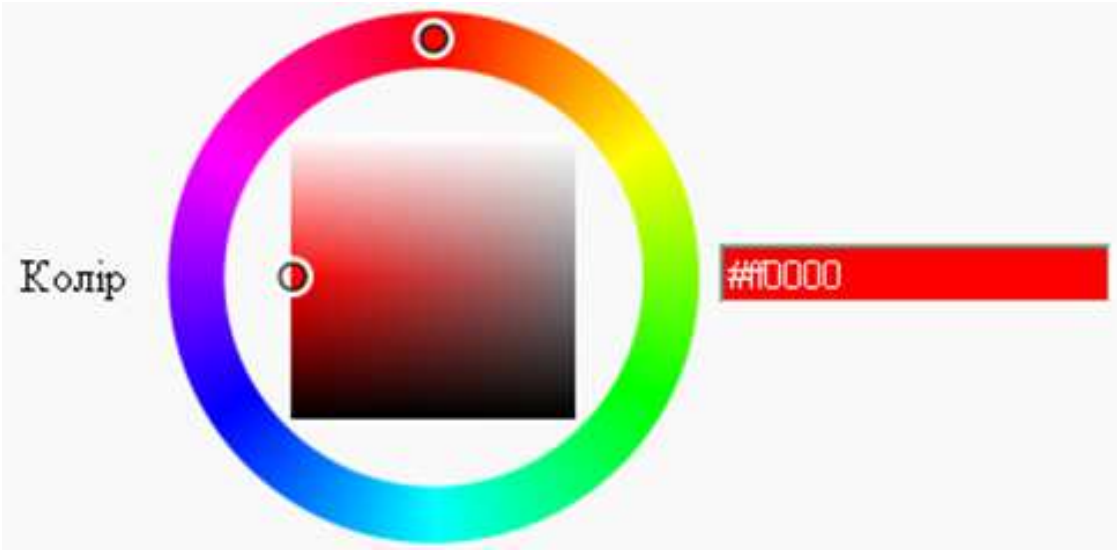

```
@interact
def _(c = color_selector ((0,1,0), "Колір", widget =
'colorpicker' ,hide_box=True)):
    show(plot(sin(7*x), color = c))
```

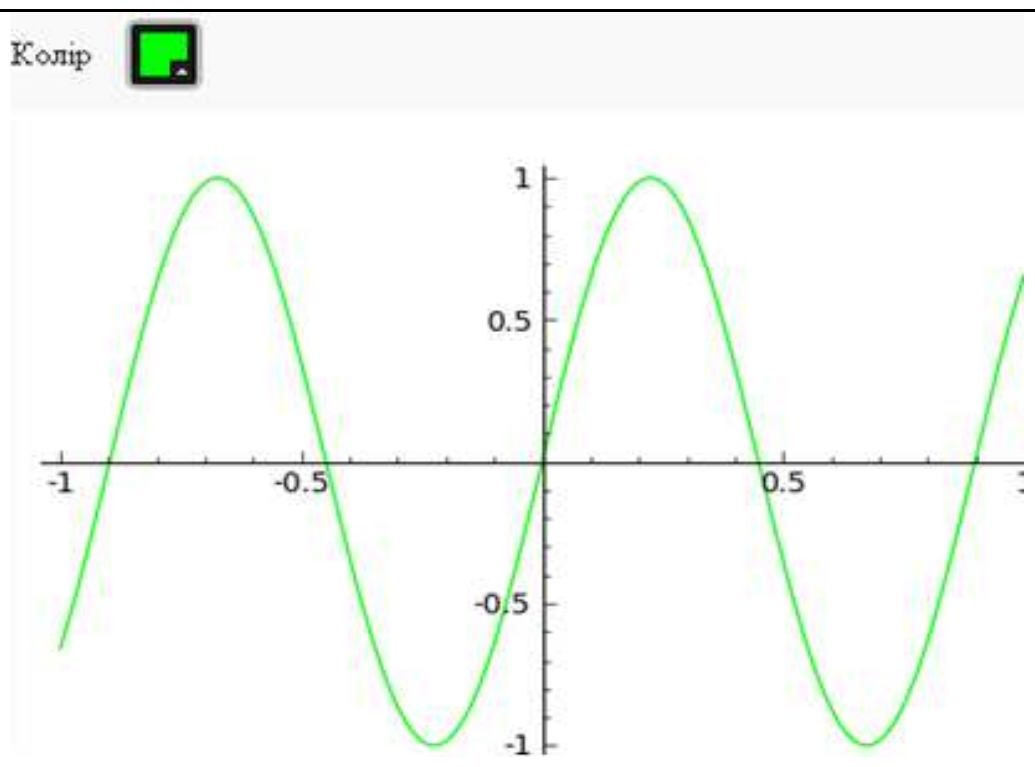

```
@interact
def _(c = color_selector((1,0,0), "Колір", widget=
'jpicker', hide_box=False)):
    show(plot(sin(7*x), color = c))
```

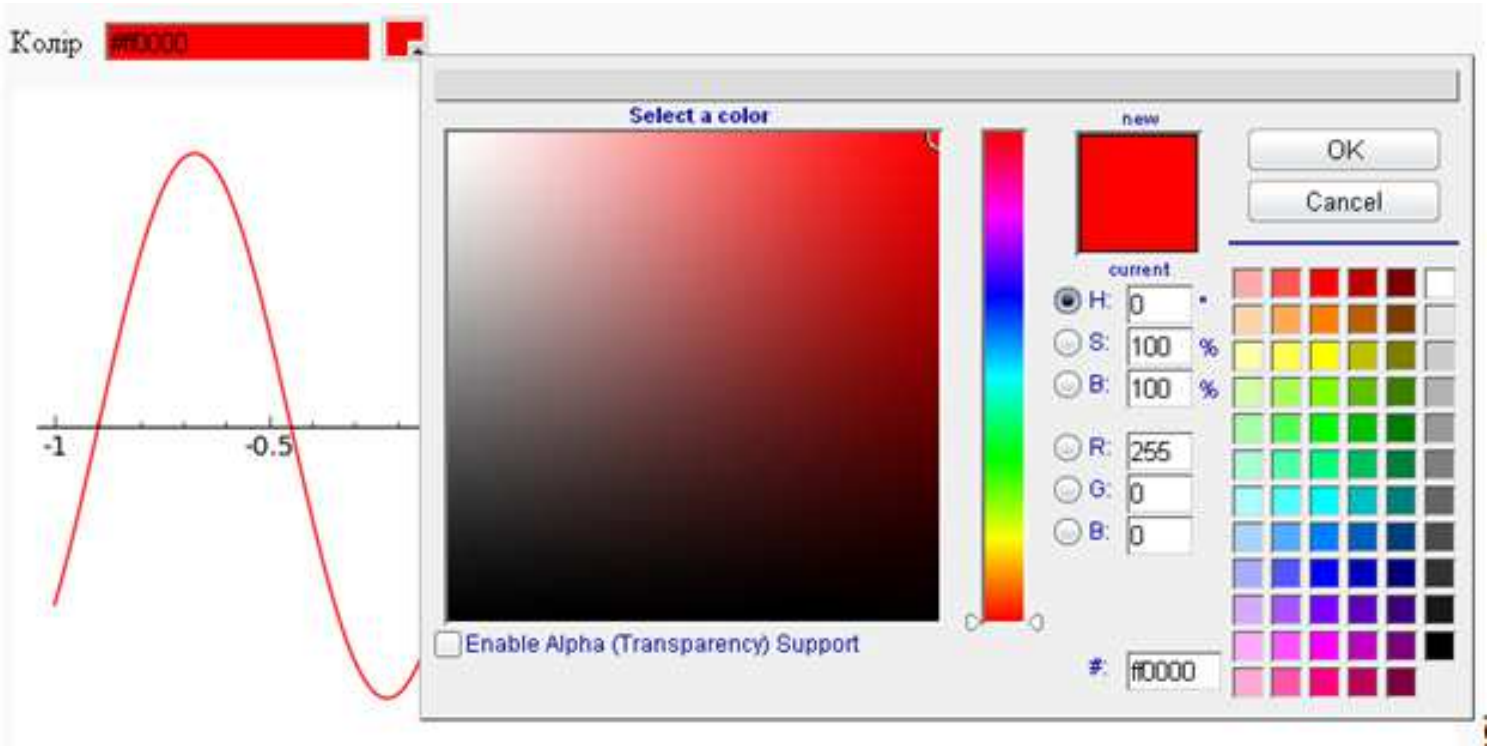

top, bottom, left, right (вгору, вниз, ліворуч, праворуч відповідно). Приклади прийомів одночасного програмування декількох елементів управління наведено у табл. 2.

*Таблиця 2*

**Приклади програмування декількох елементів управління**

**Приклад 1**

```
@interact(layout={'top':[['a','b']],'left':\
[['c']]})
def _(a=input_box("2+2", 'Вираз', width=5),
b=slider(2, 5, 3/17, 3, 'R'),
c=checkbox(False, "Відображати")):
    show(a)
```

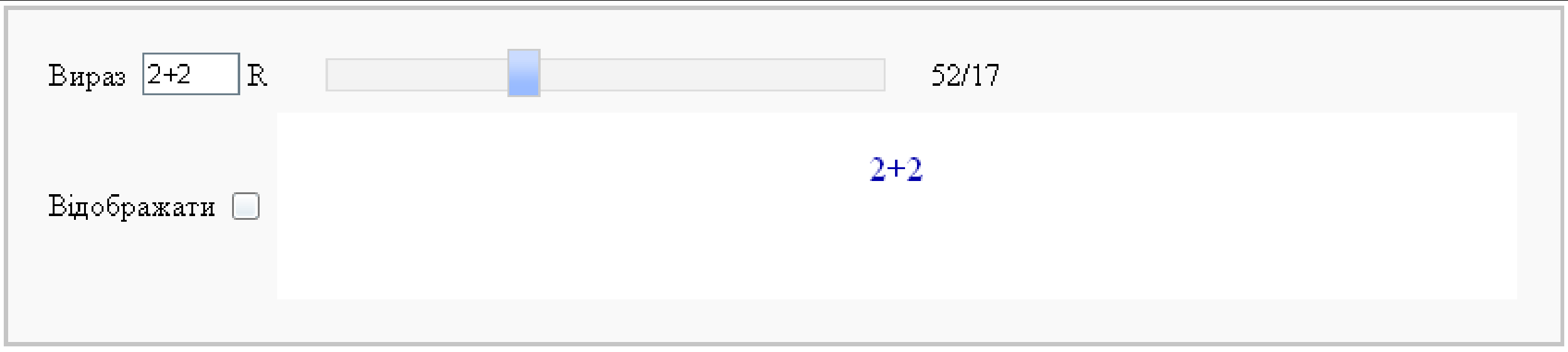

**Приклад 2.**

```
@interact(layout={'right':[['a','b']],\
'left':[['c']]})
def _(a=input_box("2+2", 'Вираз', width=5),
b=slider(2, 5, 3/17, 3, 'R'),
c=checkbox(True, "Відображати")):
    show(a)
```

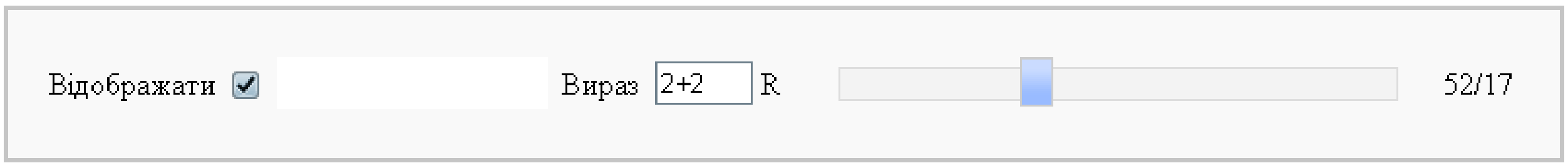

**Приклад 3.**

```
@interact(layout=[['a','b'],['c','d'],['e']])
def _(a=input_box("2+2", 'Вираз'),
b=selector([1,2,7], default=2),
c=checkbox(False, "Відображати"),
d=checkbox(True),
e=slider([1..10], None, None, 3, 'N')):
    show(a)
```

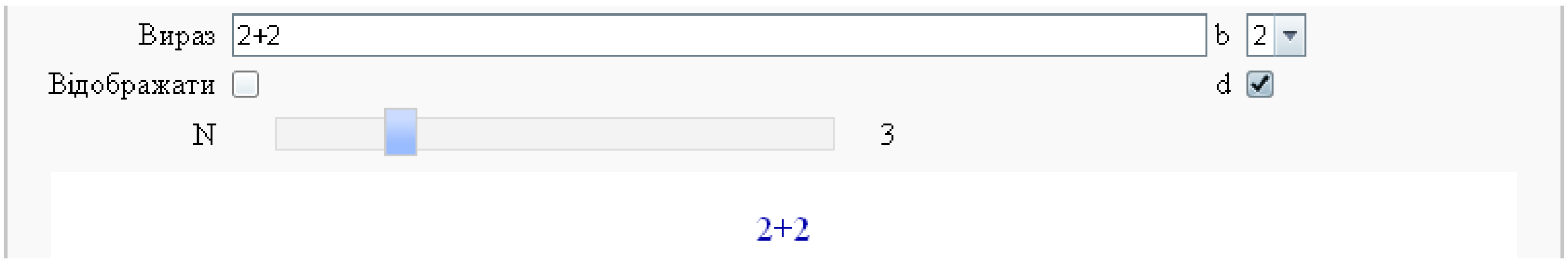

Крім елементів управління й основного програмного коду, невід'ємною частиною комп'ютерних моделей з графічним інтерфейсом і напівавтоматичним режимом управлінням є різноманітні підписи, що надають можливість більш детально пояснювати і коментувати математичні вирази. Для додавання підписів у полі графічних побудов у вигляді таблиць, кольорових графіків, текстів різного формату і підписів використовують мову HTML, а для подання виразів і формул у природній математичній нотації – мову LaTeX. Розглянемо більш детально на прикладах.

Під час оголошення функцій, що створюють елементи управління типу «Повзунок» параметру `label` (приклад 1), присвоюється значення, записане за допомогою команд HTML, що передбачають збільшення розміру і кольору тексту на екрані. Для виведення підпису «Рівняння виду …» у звичному математичному виді необхідно записати службове слово html, після якого в дужках за допомогою відповідних тегів вказати колір, тип шрифту, розмір та розташування тексту, що виводиться. Виведення функції $y = \sin(ax + b)$ здійснюється з використання команд мови LaTeX, на що вказують одинарні лапки й обмежувальними знаки «$». Кожна команда розпочинається символом «\» (backslash – обернений слеш), після якого зазначається власне ім'я команди: команда «\sin» відображає функцію *sin*, команда «\cdot» – операцію множення.

**Приклад 1.**

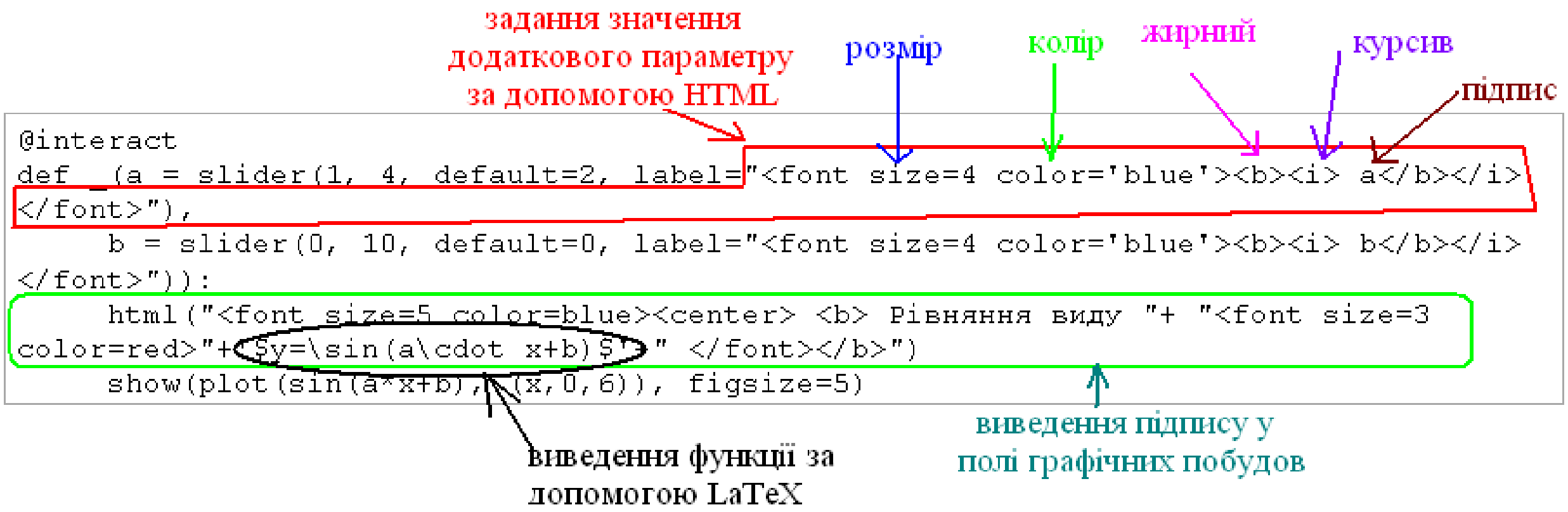

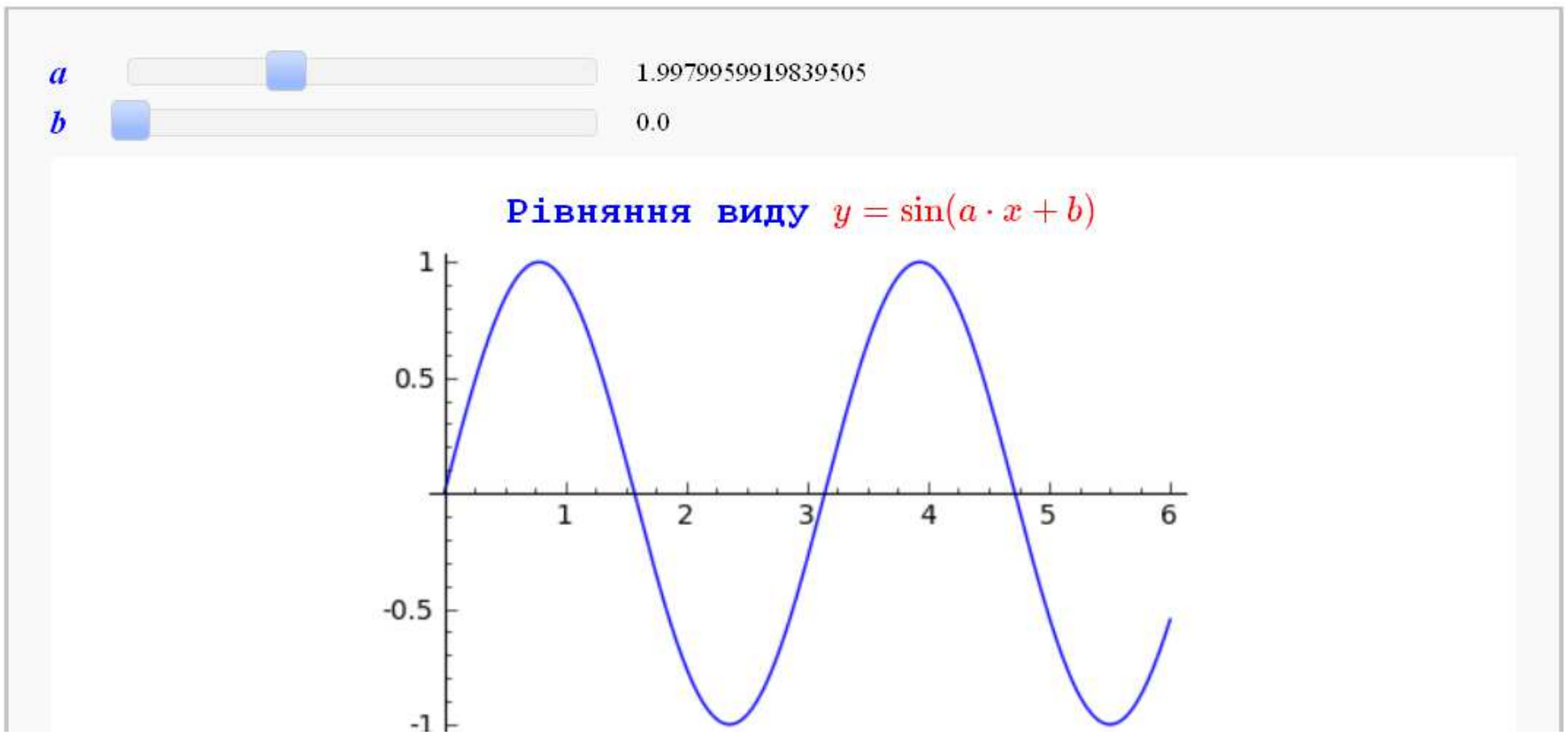

На відміну від прикладу 1, приклад 2 передбачає разом із зміною параметра в полі «Введіть функцію» зміну функції, що виводиться разом з підписом «Зображено графік функції:». Це досягається за рахунок використання такої комбінації:

+'$f(x)=%s$'%latex(f)+, де f(x) – назва функції (постійне значення), %s – формат для виведення рядкової змінної, latex(f) – функція перетворення f у природну математичну нотацію.

**Приклад 2.**

функція для створення елемента управління «Поле для введення»

параметр для задання надпису ліворуч від елемента

основний параметр для задання значення, що повертається функцією за замовчуванням

```
@interact
def _(f=input_box(label="Введіть функцію:",default=sin(x)*e^(-x))):
    var('x')
    show(plot(f,0,15,rgbcolor='green', thickness=1))
    html("<b><font color=blue size=4>Зображено графік функції:
"+'$f(x)=%s$'%latex(f)+" "+"</font></b>")
```

виведення підпису у полі графічних побудов

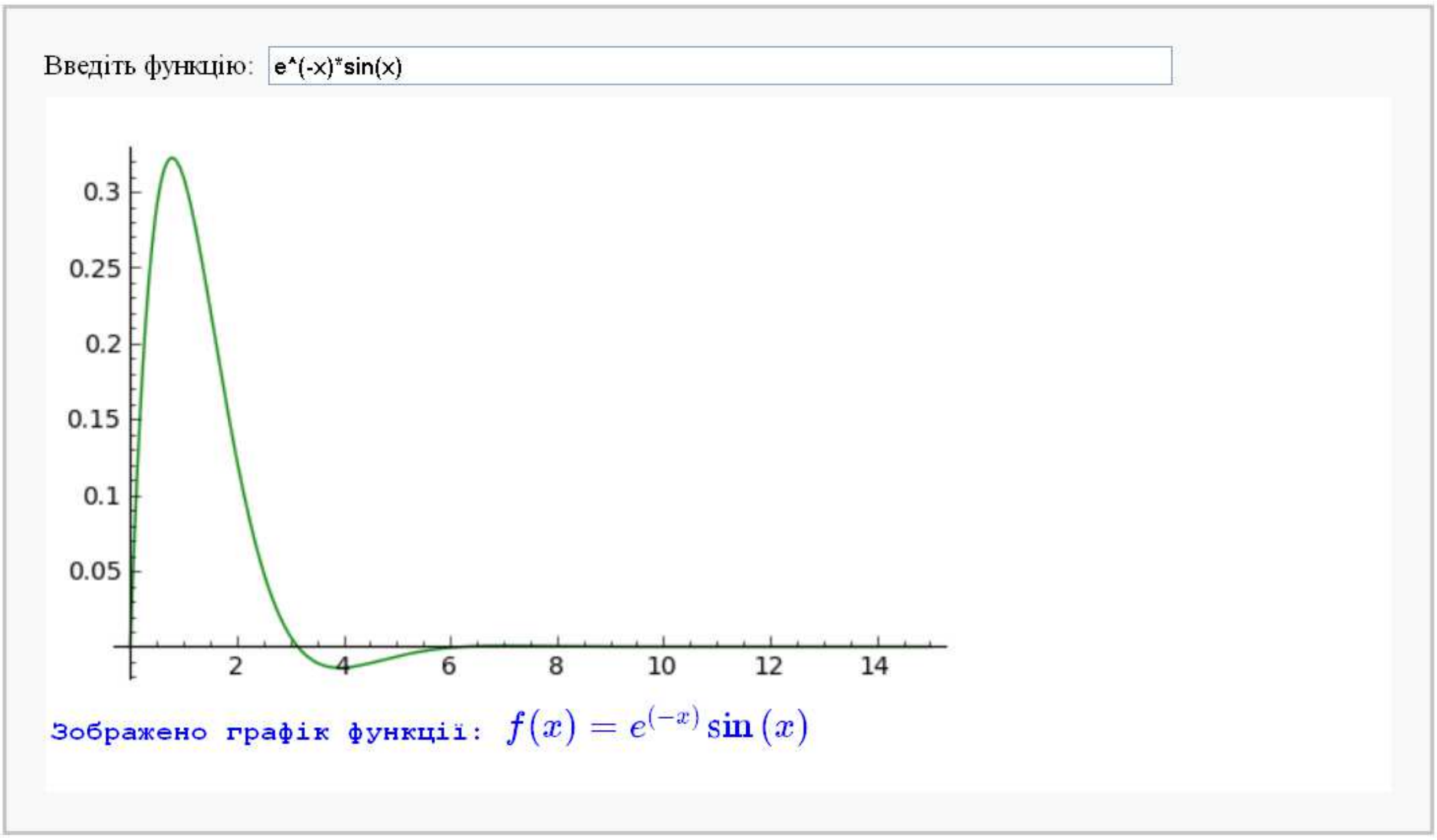

Отже, для створення комп’ютерних моделей з графічним інтерфейсом і напівавтоматичним режимом управлінням необхідно:

– відповідно до моделі, що розробляється, створити елементи управління;

– описати програмний код моделі, що створюється, за допомогою мови програмування Python;

– додати підписи різного формату у вигляді тексту або математичних виразів за допомогою команд мови HTML і LaTeX.

Зазначимо, що перед початком роботи у потрібній комірці робочого аркуша,

можна задати такі модифікатори : `%hide` – приховувати вміст комірки до першого звернення до неї; `%auto` – автоматично виконувати вміст комірки під час відкриття аркуша; `%html`, `%sage`, `%python`, `%maxima` – використовувати обраний інтерпретатор для опрацювання вмісту комірки (за замовчуванням `%sage`).

Продемонструємо методику побудови комп’ютерних моделей з графічним інтерфейсом і напівавтоматичним режимом управлінням на прикладі створення генератора завдань контрольної роботи з лінійної алгебри (рис. 1). Актуальність створення програм-генераторів навчальних завдань має підтвердження у роботі [8], у якій автори описали методику створення банку задач за допомогою СКМ Maple, що надає викладачеві гнучкий інструмент для формування різноманітних завдань згідно сучасних вимог щодо здійснення контролю у процесі вивчення вищої математики.

Кількість варіантів 2

**Варіант 1**

1. Обчисліть визначник $\begin{vmatrix} -1 & 1 & -2 \\ -2 & 2 & -1 \\ -1 & 3 & 4 \end{vmatrix}$, використовуючи: а) правило трикутників; б) метод розкладання визначника за елементами деякого рядка або стовпця; в) метод зведення до трикутного вигляду.

2. Обчисліть добуток **AB** та **BA**:
**A**=$\begin{pmatrix} -4 & -4 \end{pmatrix}$, **B**=$\begin{pmatrix} 0 & 4 & 3 \\ 3 & 2 & 1 \end{pmatrix}$

3. Знайдіть матрицю, обернену до $\begin{pmatrix} -4 & -2 & -1 \\ 2 & -3 & -4 \\ 3 & 4 & -3 \end{pmatrix}$:

а) за допомогою алгебраїчних доповнень;

б) методом Гауса.

4. Розв'яжіть матричне рівняння: **ABX=C**:
**A**=$\begin{pmatrix} 5 & -5 \\ -2 & 4 \end{pmatrix}$, **B**=$\begin{pmatrix} 2 & -2 \\ -3 & -4 \end{pmatrix}$, **C**=$\begin{pmatrix} -4 & 3 \\ 1 & 0 \end{pmatrix}$

5. Обчисліть $f(\mathbf{A})$, якщо **A**=$\begin{pmatrix} 4 & 2 & 1 \\ 4 & -2 & -1 \\ 0 & -5 & -4 \end{pmatrix}$, $f(x) = -3x^2 + 5x - 5$

6. Розв'яжіть систему рівнянь: а) методом Крамера; б) методом оберненої матриці; в) методом Гауса.
$$\begin{cases} 3x_1 - 2x_2 - 5x_3 = -24 \\ 4x_1 - 4x_2 - 3x_3 = -11 \\ -5x_1 - 4x_2 = 13 \end{cases}$$

*Рис. 1. Інтерфейс генератора контрольної роботи з лінійної алгебри*

Відповідно до рис. 1 даний генератор контрольної роботи містить поле для вибору кількості варіантів, що використовується як лічильник циклів, тому за допомогою функції `input_box` створимо відповідний елемент управління. Для можливості швидкої перевірки контрольної роботи, під час створення генератора завдань необхідно відразу оголосити порожній масив відповідей (`answers`), до якого будемо записувати відповіді всіх завдань контрольної роботи. Потім організуємо цикл за кількістю варіантів, у якому будемо виводити повідомлення про номер варіанту і його завдання:

```
    @interact
    def   mkr1(numvar   =   input_box(default=2,   label=
"Кількість варіантів")):
        answers=[]
        for var in range(numvar):
            html("</pre>  <center>  <font  size=+2>  Варіант
%s</font></center>"%(var+1))
```

створення поля для введення кількості варіантів (за замовчуванням - 2)

масив відповідей

цикл за кількістю варіантів

виведення повідомлення про номер варіанту

Для генерації *першого* завдання контрольної необхідно задати матрицю розмірності 3 на 3 і заповнити її цілими випадковими числами (у даному випадку від –5 до 5). Водночас відразу обчислюється визначник матриці, значення якого записується у масив відповідей.

```
A=matrix(QQ,3,3)
for i in range(3):
    for j in range(3):
        A[i,j]=randint(-5,5)
answers.append(det(A))
```

Слід зазначити, що функція `show()` не надає можливості виведення визначника у математичній нотації (матриці з вертикальними дужками), тому виникла необхідність створення допоміжної функції (`nicedet`), що за матрицею `M` формує послідовність команд мовою LaTeX, яка зображує матрицю з вертикальними дужками.

```
def nicedet(M):
    s="$\\left|\\begin{array}{ccc} "
    for i in     range(M.nrows()):
        for j in range(M.ncols()):
            s=s+" "+latex(M[i,j])
```

```
            if j!=M.ncols()-1:
               s=s+"&"
        s=s+"\\\\"
    s=s+"\\end{array}\\right|$"
    return s
```

Як параметр даної функції передається матриця `M`, результат перетворення якої записують у рядкову змінну `s`.

Для зображення самого визначника і завдань, які необхідно виконати, потрібно скористатися тегами HTML:

```
html("<p>1. Обчисліть визначник %s, використовуючи: "%nicedet(A))
html("а) правило трикутників; ")
html("б) метод розкладання визначника за елементами деякого рядка або стовпця; ")
html("в) метод зведення до трикутного вигляду.<br>")
```

Для генерації *другого* завдання випадково задаємо дві матриці наперед невідомої розмірності (кількість рядків і стовпців у діапазоні від 1 до 3) і заповнюємо їх випадковими цілими числами від –5 до 5. Так як розмірність матриць наперед невідома, то операція множення не завжди є виконуваною, тому скористаємося блоком «`try ... except`». Отже, якщо добуток матриць існує, то його значення записуємо у масив відповідей, а якщо матриці перемножити неможливо, то у масив відповідей записуємо рядкову змінну з відповідним повідомленням.

```
try:
  res=A*B
except:
  res="Добуток AB не існує"
answers.append(res)
try:
  res=B*A
except:
  res="Добуток BA не існує"
answers.append(res)
```

Отже, у результаті виконання операції множення матриць у масив відповідей запишеться або дві матриці, або два повідомлення, або матриця і повідомлення. За аналогічним принципом генерують *третє*, *четверте* та *п'яте* завдання.

У *шостому* завданні необхідно розв'язати систему лінійних рівнянь трьома способами. Для цього спочатку задаємо вектор із 3-х цілих випадкових чисел і одразу додаємо його до масиву відповідей. Отже, замість того, щоб шукати розв'язки системи, спочатку задамо її корені, і за ними будуємо систему рівнянь. Це необхідно для того, щоб у процесі виконання завдання кожна утворена система лінійних рівнянь завжди мала розв'язок, при цьому корені системи були цілими числами. Потім формуємо матриці `A` і `B`, зокрема у матрицю `B` записуємо числа з матриці `A` як коефіцієнти і множимо їх на відповідні корені. Отже, `A` – коефіцієнти лівої частини системи, а `B` – правої. Далі формуємо вектор із трьох символьних значень `x1, x2, x3` (змінних) і матрицю-стовпець `X` – вектор невідомих. Далі, за допомогою команд LaTeX формуємо рядкову змінну і виводимо матриці `A*X` і `B` поелементно таким способом, щоб зліва була одна велика фігурна дужка, а в кожному рядку стояв знак «=».

```
res=[randint(-5,5),randint(-5,5),randint(-5,5)]
k=vector(res)
answers.append(k)
A=matrix(QQ,3,3)
for i in range(A.nrows()):
    for j in range(A.ncols()):
        A[i,j]=randint(-5,5)
B=matrix(QQ,3,1)
for i in range(A.nrows()):
    for j in range(A.ncols()):
        B[i,0]=B[i,0]+A[i,j]*res[j]
L=[x1,x2,x3]
X=matrix(3,1,L)
s="$\\left \\{ \\begin{eqnarray}"
for i in range(3):
    s=s+"%s=%s\\\\" %(latex((A*X)[i][0]), latex(B[i][0]))
```

```
s=s+"\\end{eqnarray} \\right.$"
html(s)
```

Для відокремлення завдань від відповідей горизонтальною рискою використовуємо відповідний тег HTML. Потім організуємо цикл для виведення всіх відповідей кожного варіанту з масиву `answers`:

```
html("<hr>")          # лінія, що розділяє відповіді
html("<p><b>ВІДПОВІДІ</P></b>")
html("<p><p><p><p><p><p><p><p>")
for var in range(numvar):
  html("<p><b>Варіант %s</b>: "%(var+1))
  html("1.$%s$;2.<b>AB</b>=$%s$,<b>BA</b>=$%s$;3.   $%s$;
4. $%s$; 5. $%s$; 6. $%s$"% \
         (latex(answers[var*7+0]),\
          latex(answers[var*7+1]),
latex(answers[var*7+2]),latex(answers[var*7+3]),
latex(answers[var*7+4]), \
          latex(answers[var*7+5]), latex(answers[var*7+6])))
```

Отже, за допомогою даного генератора викладачеві надається можливість не тільки отримати необхідну кількість варіантів завдань, а й відповіді до них.

Наведений вище приклад генератора навчальних завдань обов'язково виводить відповіді на екран після всіх варіантів завдань, проте, за бажання, за допомогою елемента управління «прапорець», можна вибрати, показувати відповіді, чи ні (рис. 2).

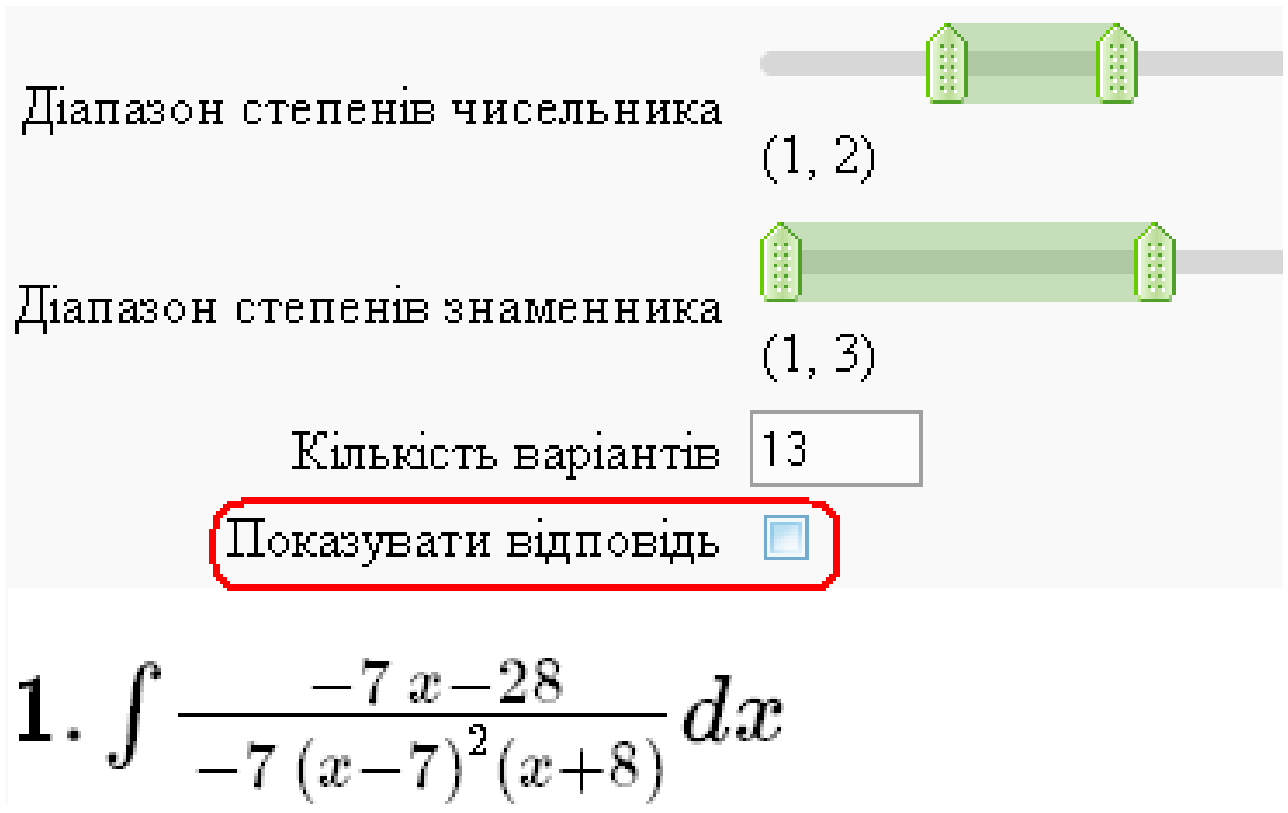

*Рис. 2. Інтерфейс генератора завдань з теми*

*«Інтегрування раціональних функцій»*

Одним із компонентів ММС «Вища математика» є навчальні експертні системи, створення бази знань яких відбувається за допомогою eXpertise2Go.

Експертна система eXpertise2Go [9; 10; 11] є вільно поширюваним Web-орієнтованим програмним засобом (Web-НЕС), що надає можливість генерувати базу знань у форматі e2gRuleEngine за допомогою інструменту для створення і перевірки таблиць розв’язків e2gRuleWriter.

Для створення бази знань потрібно заповнити таблицю розв’язків через виконання таких дій:

1) завантажити файл e2gRuleWriter.jar, у результаті чого на екрані з’явиться таблиця (рис. 3), що зображує три режими («**CONDITIONS**» – умови, «**ACTIONS**» – дії, «**Rule**» – правила) різними кольорами: жовтий, зелений та голубий для «CONDITIONS», «ACTIONS» та «Rule» відповідно, синій колір – ознака активності рядка чи стовпця;

2) за допомогою опцій поля «**Prompt Type**» визначити, як дані будуть запитуватися у користувача: *YesNo* – вибір типу «так–ні», *MultChoice* – створює декілька рядкових вхідних повідомлень і додає повідомлення типу «не можу відповісти» (у вигляді вибору кнопок), *ForcedChoice* – створює декілька рядкових вхідних повідомлень, *Choice* – вибір, подібний до *MultChoice*, тільки варіанти відповідей задаються за допомогою списку, що розкривається, *AllChoice* – використовують коли необхідно вибрати більше ніж один варіант відповіді, *Numeric* – для числового введення відповіді;

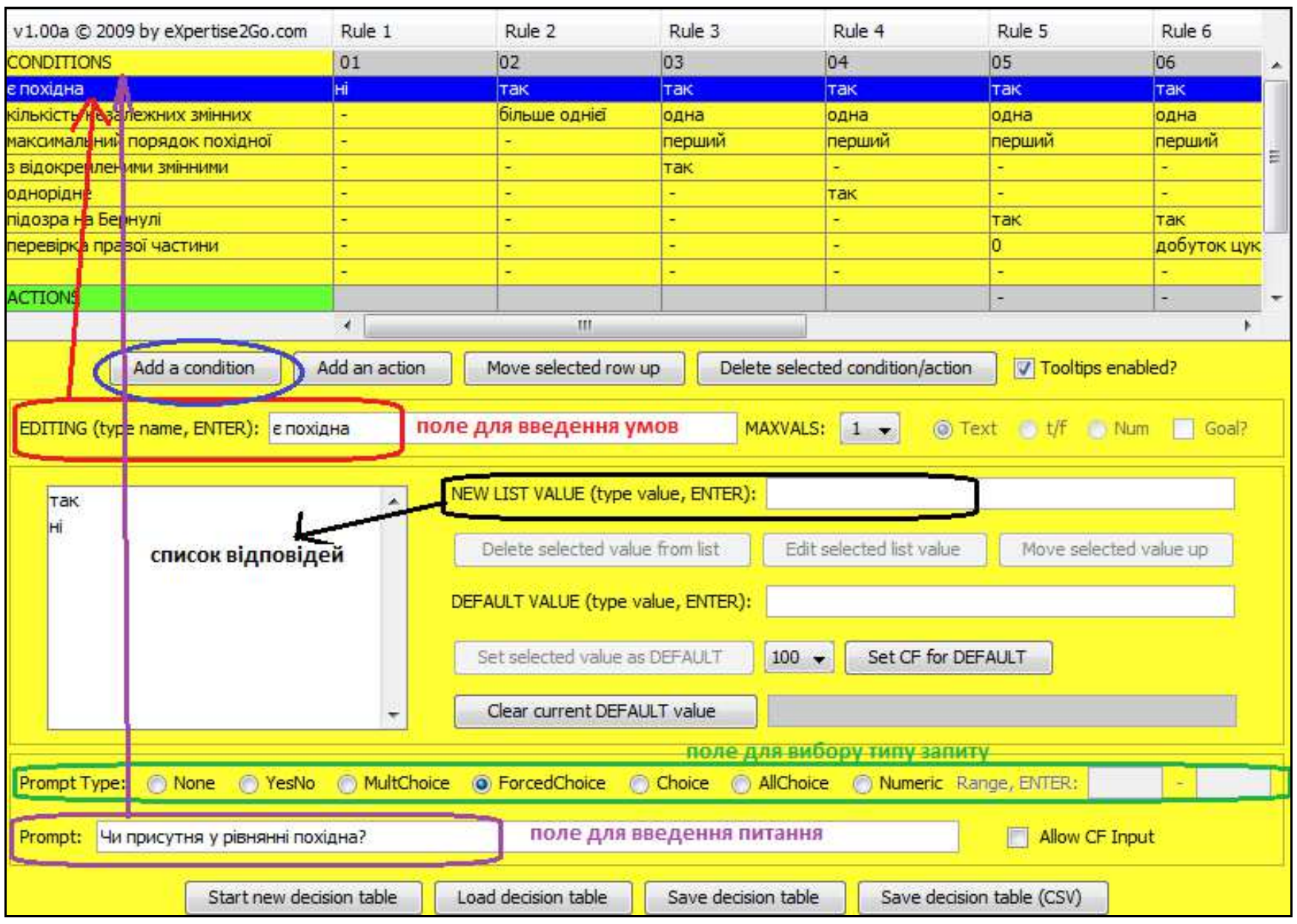

*Рис. 3. Інтерфейс таблиці для визначення «CONDITIONS»*

3) визначити умови і запитання: для цього натиснути лівою кнопкою миші на рядку поля «CONDITIONS» (воно стане синього кольору), після чого у полі «EDITING» з’явиться курсор, що надає можливість увести текст умови. Для відображення необхідного запитання (відповідно до введеної умови) в інтерфейсі eXpertise2Go, потрібно у полі «**Prompt**» увести текст запитання, що буде відображатися користувачеві. Для введення варіантів відповідей потрібно скористатися полем «**NEW LIST VALUE**», при цьому список всіх відповідей буде відображатися у відповідному вікні. У випадку, коли необхідно додати ще кілька умов, а рядків не вистачає, потрібно скористатися кнопкою *Add a condition*;

4) визначити дії (вказати висновок, рекомендацію): для цього натиснути лівою кнопкою миші на рядку поля «ACTIONS» (воно стане синього кольору), після чого у полі «EDITING» з’явиться курсор, що надає можливість увести частину тексту загального для всіх висновків (у даному випадку це «тип рівняння»). Щоб встановлення типу рівняння було «метою» консультації даної експертної системи, необхідно поставити прапорець у полі «**Goal**». Для введення варіантів висновків потрібно скористатися полем «NEW LIST VALUE», при цьому список всіх висновків (типів рівнянь) буде відображатися у відповідному вікні (рис. 4);

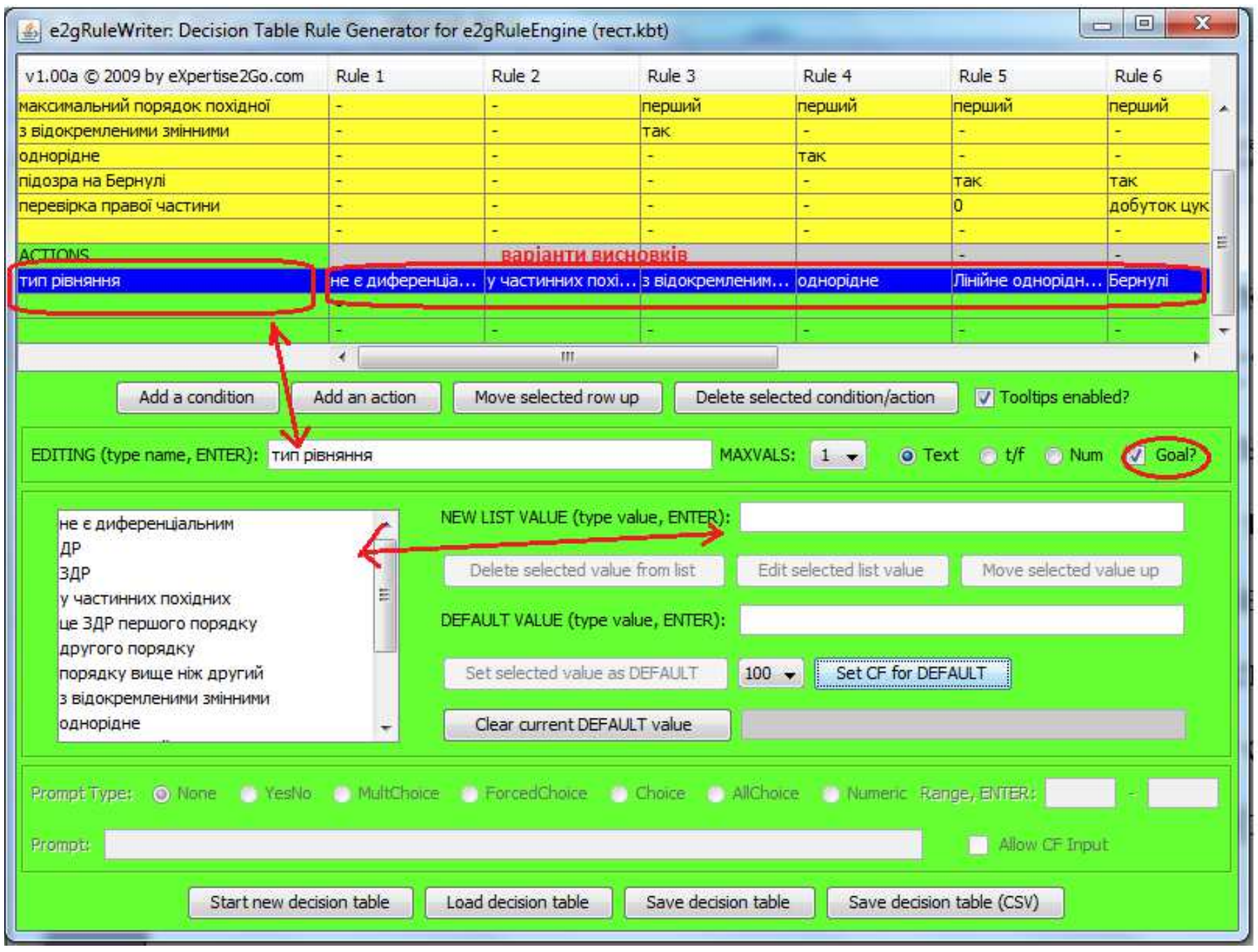

*Рис. 4. Інтерфейс таблиці для визначення дій НЕС*

5) визначити правила за допомогою полів «Rule». Спочатку потрібно обрати необхідне правило (наприклад Rule 5) і полі «EDITING» вказати його ім’я (у даному випадку – 05). Потім у кожному рядку правила відповідно до кожної умови із списку, що випадає, необхідно вибрати варіант відповіді (натиснути мишею на відповідний рядок). Аналогічні дії слід виконати для поля «ACTIONS». За необхідності вибрати опції для спрощення і перевірки визначених правил (рис. 5);

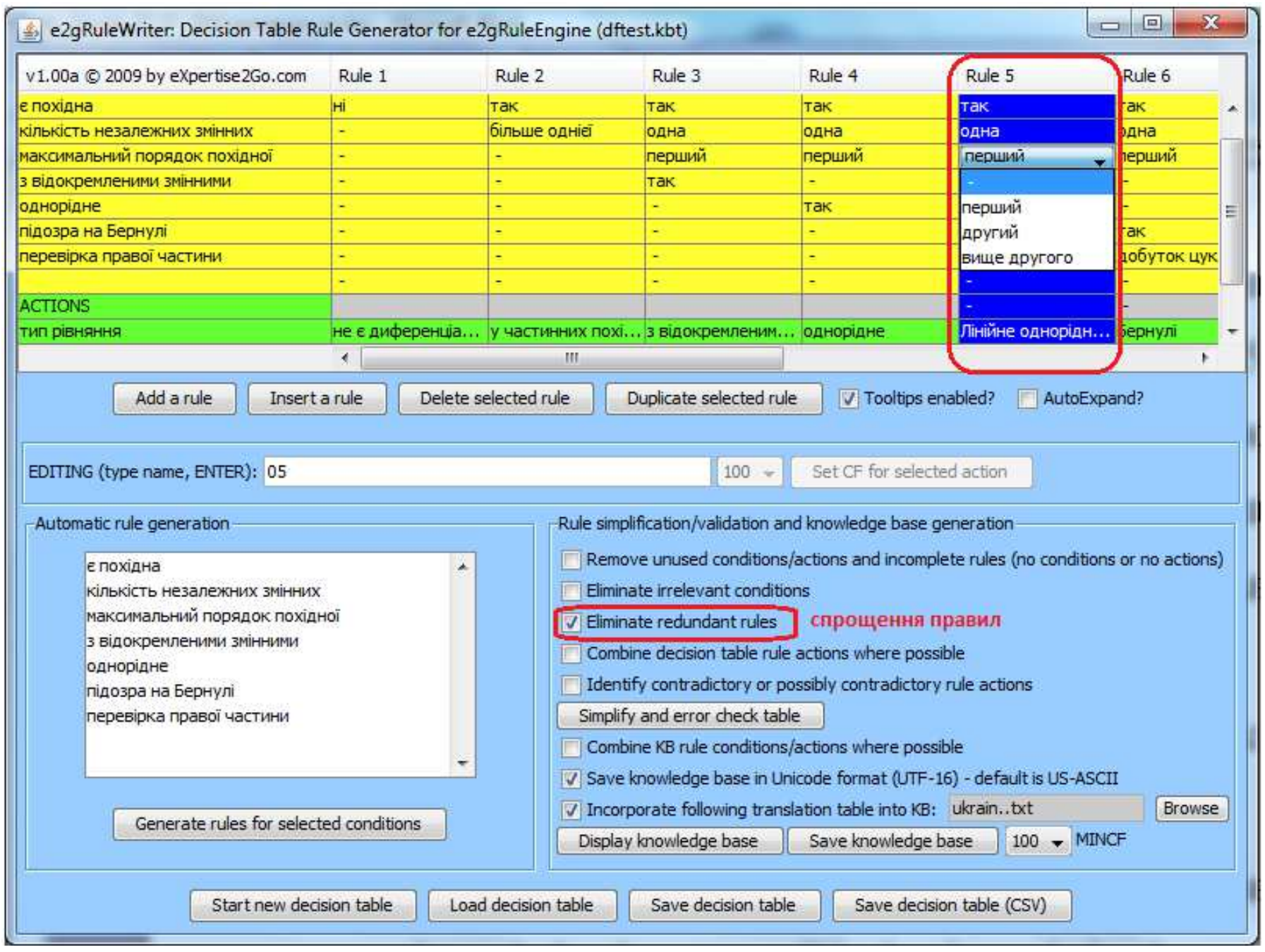

*Рис. 5. Інтерфейс таблиці для визначення правил НЕС*

6) згенерувати і переглянути базу знань натисканням кнопок «Save knowledge base» і «Display knowledge base» відповідно (рис. 6). Також, на цьому етапі слід обрати фактор вірогідності (поле «**MINCE**»), тобто з яким відсотком довіри експертна система надає рекомендацію. Зауважимо, що оскільки правила, дії та умови написані українською мовою, то збереження бази знань експертної системи необхідно виконати у форматі Unicode. Для відображення у головному вікні експертної системи кнопок управління українською мовою, виконано локалізацію, тому під час збереження бази знань необхідно завантажити відповідний файл у форматі «txt»:

```
REM Надписи на кнопках
TRANSLATE B_SUBMIT = "Відповісти"
TRANSLATE B_EXPLAIN = "Пояснити"
TRANSLATE B_WHYASK = "Чому питаємо?"
TRANSLATE B_RESTART = "До початку"
TRANSLATE B_RETURN = "Повернутися"
REM Повідомлення
TRANSLATE TR_KB = "База знань:"
TRANSLATE TR_NORESP = "Не знаю"
```

```
TRANSLATE TR_HOWCONF = "Наскільки Ви впевнені у відповіді?"
TRANSLATE TR_LOWCONF = "Наполовину (50%)"
TRANSLATE TR_HICONF = "Цілком (100%)"
TRANSLATE TR_YES = "Так"
TRANSLATE TR_NO = "Ні"
REM TRANSLATE TR_FALSE = "хиба"
TRANSLATE TR_RESULTS = "ВИСНОВОК:"
TRANSLATE TR_MINCF = "Мінімальний коефіцієнт довіри для прийняття значення як факту:"
TRANSLATE TR_NOTDETERMINED = "неможливо визначити"
TRANSLATE TR_ISRESULT = "є:"
TRANSLATE TR_WITH = "з"
TRANSLATE TR_CONF = "% довіри"
TRANSLATE TR_ALLGOALS = "всі висновки"
TRANSLATE TR_VALUE = "Значення"
TRANSLATE TR_OF = ""
TRANSLATE TR_THISRULE = "Відповідь для цього правила була уведена з коефіцієнтом довіри "
TRANSLATE TR_RULEASGN = "і надано значення"
TRANSLATE TR_TOFIND = "Для знаходження"
TRANSLATE TR_AVALUE = "значення для"
TRANSLATE TR_ISNEEDED = "необхідно випробувати дане правило:"
TRANSLATE TR_RULE = "ПРАВИЛО:"
TRANSLATE TR_IF = "Якщо"
TRANSLATE TR_THEN = "То"
TRANSLATE TR_AND = "і"
TRANSLATE TR_OR = "або"
TRANSLATE TR_EQUAL = "-"
TRANSLATE TR_LESSTHAN = "менше, ніж"
TRANSLATE TR_GREATER = "більше, ніж"
TRANSLATE TR_NOTEQUAL = "не дорівнює"
```

```
TRANSLATE TR_VALUEFOR = "Значення для:"
TRANSLATE TR_FOUND = "було визначено"
TRANSLATE TR_NOTFOUND = "не було визначено"
TRANSLATE TR_WASINPUT = "було уведено з "
TRANSLATE TR_DETERMINED = "Визначено"
TRANSLATE TR_IS = "-"
TRANSLATE TR_FROM = "з:"
TRANSLATE TR_DEFAULTED = "було встановлено за замовчуванням у"
TRANSLATE TR_ONE = "одне зі значень"
TRANSLATE TR_HOWCF1 = "Обчислення % довіри за кількома джерелами для"
```

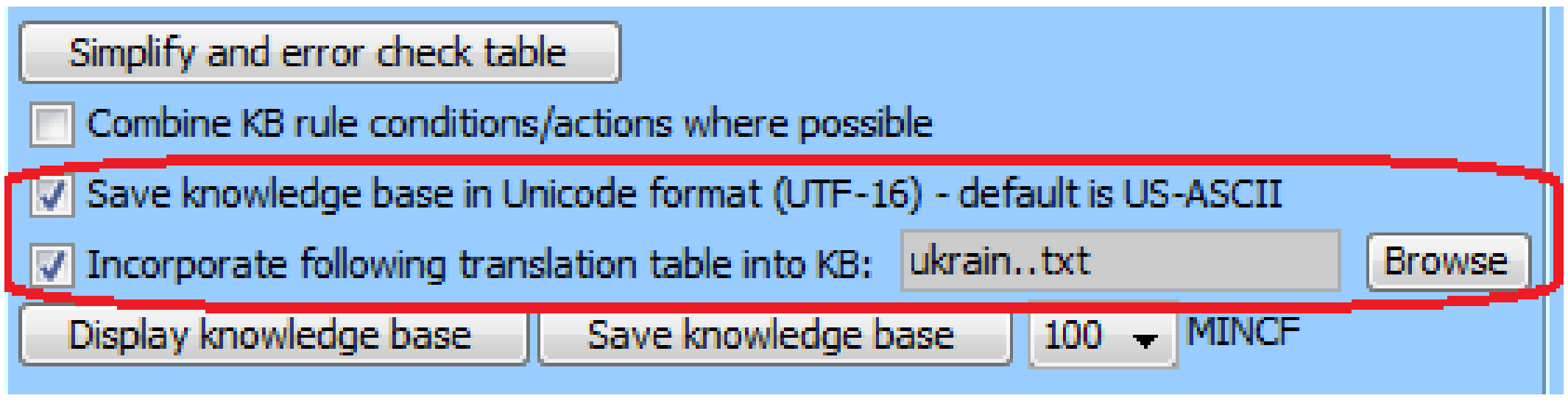

*Рис. 6. Дії, необхідні для генерації бази знань у форматі e2gRuleEngine*

Зазначимо, що після натискання кнопки «Save knowledge base» отримаємо файл з розширенням «kb», який потрібно зберегти у кодуванні UTF-16;

7) випробувати базу знань за допомогою аплету e2gRuleEngine.jar. Для цього у треба створити веб-сторінку за допомогою наступних тегів і зберегти файл з розширенням «html»:

```
<html>
<head><title>Експертна система</title></head>
<body bgcolor="#c0c0c0">
<center>
<applet code="e2gRuleEngine.class" archive="e2gRuleEngine.jar" width=800 height=500>
<param name="KBURL" value="dftest2.kb">
<param name="APPTITLE" value="Визначення типу диференціального рівняння">
```

```
<param name="APPSUBTITLE" value="К. І. Словак">
<param name="BGCOLOR" value="#ffff00">
<param name="TITLECOLOR" value="#0000aa">
<param name="PROMPTCOLOR" value="#0000ff">
<PARAM NAME="PROMPTSIZE" VALUE="15">
<param name="WORDWRAP" value="true">
<param name="KBENCODING" VALUE="UTF-16">
<param name="DEBUG" value="false">
<PARAM NAME="MONOFONT" VALUE="Arial Unicode MS">
<PARAM NAME="LOADMSG" VALUE="Завантаження бази знань...">
<PARAM NAME="NOLOGO" VALUE="TRUE">
<param name="STARTBUTTON" value="Увійти до ЕС">
</applet>
</body>
</html>
```

Отже, для перевірки роботи створеної експертної системи потрібно завантажити html файл, який повинен бути розташований в одному каталозі з файлом бази знань та аплетом e2gRuleEngine.jar (за наявності цих трьох файлів у одному каталозі, експертна система буде працювати будь-де).

Вище було розглянуто основні етапи побудови і перевірки бази знань експертної системи eXpertise2Go, що може функціонувати на будь-якому персональному комп’ютері. У той же час, для налаштування роботи експертної системи у ММС необхідно (рис. 7):

1) створити новий робочий аркуш;

2) у активній комірці вибрати інтерпретатор html і ввести програмний код наведений вище без тегів, що відповідають за назву вікна Web-сторінки:

```
%auto   # автоматичне виконання при завантаженні
%hide   # не показувати програмний код
%html
<center>
<applet code="e2gRuleEngine.class" ar-
chive="e2gRuleEngine.jar" width=800 height=500>
<param name="KBURL" value="dftest2.kb">
```

```
<param name="APPTITLE" value="Визначення типу диференціального рівняння">
<param name="APPSUBTITLE" value="К. І. Словак">
<param name="BGCOLOR" value="#ffff00">
<param name="TITLECOLOR" value="#0000aa">
<param name="PROMPTCOLOR" value="#0000ff">
<PARAM NAME="PROMPTSIZE" VALUE="15">
<param name="WORDWRAP" value="true">
<param name="KBENCODING" VALUE="UTF-16">
<param name="DEBUG" value="false">
<PARAM NAME="MONOFONT" VALUE="Arial Unicode MS">
<PARAM NAME="LOADMSG" VALUE="Завантаження бази знань...">
<PARAM NAME="NOLOGO" VALUE="TRUE">
<param name="STARTBUTTON" value="Увійти до ЕС">
</applet>
```

3) у список «Данные» завантажити файл бази знань та аплет e2gRuleEngine.jar.

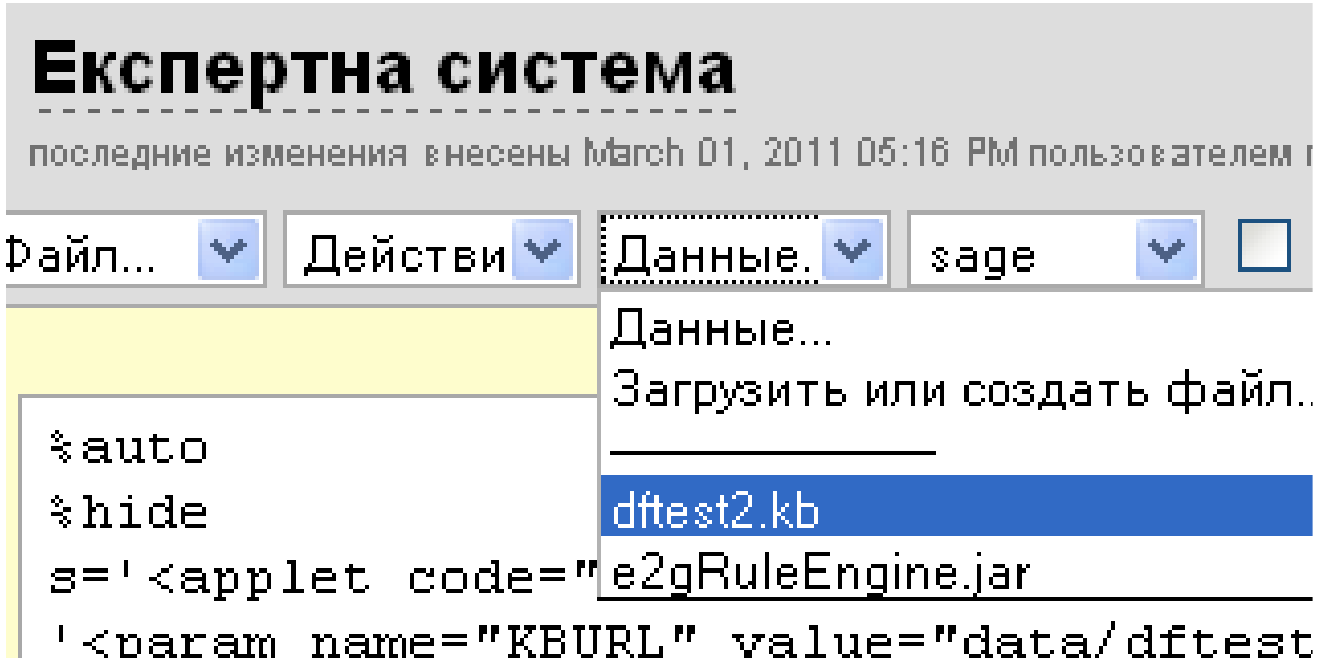

а)

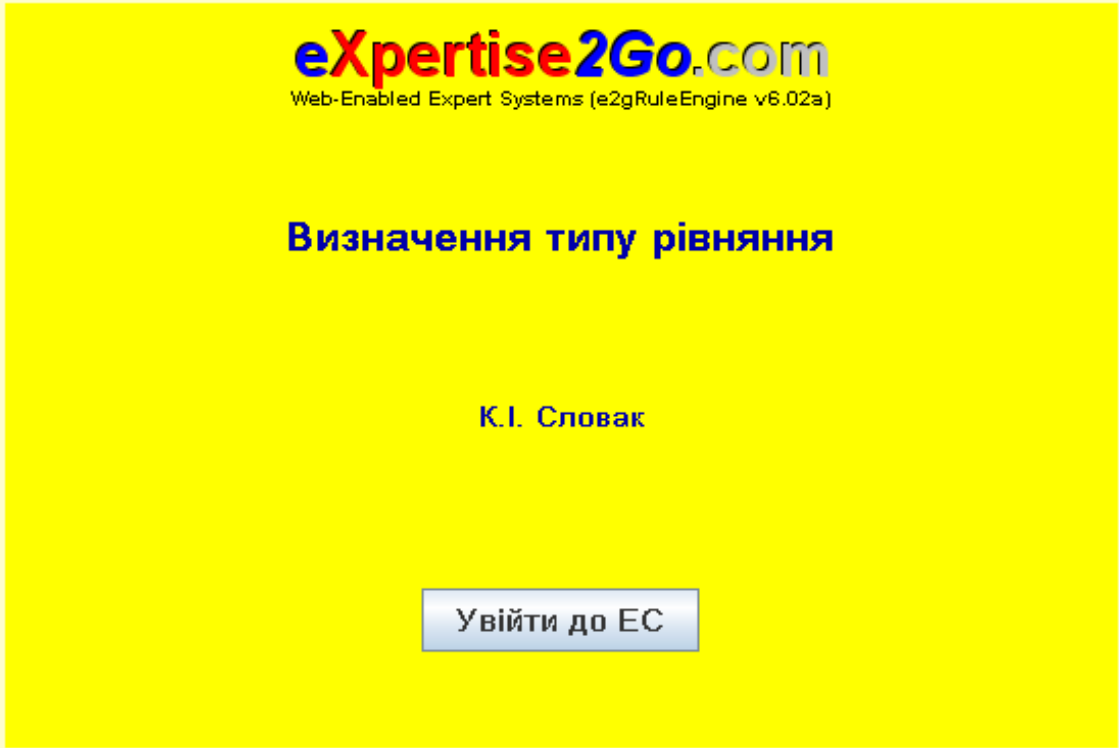

б)

*Рис. 7. Налаштування роботи експертної системи у ММС*

Створення ІДЗ, прикладів розв’язування завдань, задач для практичних занять тощо передбачає використання текстового редактора TinyMCE, убудованого в Sage. Кожен файл – робочий аркуш Sage, що містить текстові дані (навчальні завдання і приклади їх розв’язування) і прямокутні комірки з набором команд для виконання відповідних математичних розрахунків. Для відкриття TinyMCE необхідно навести мишу над коміркою, щоб з’явилася фіолетова жирна горизонтальна лінія. Потім натиснути `<Shift>` і ліву кнопку миші, у результаті чого з’явиться вікно редагування тексту, зображене на рис. 8.

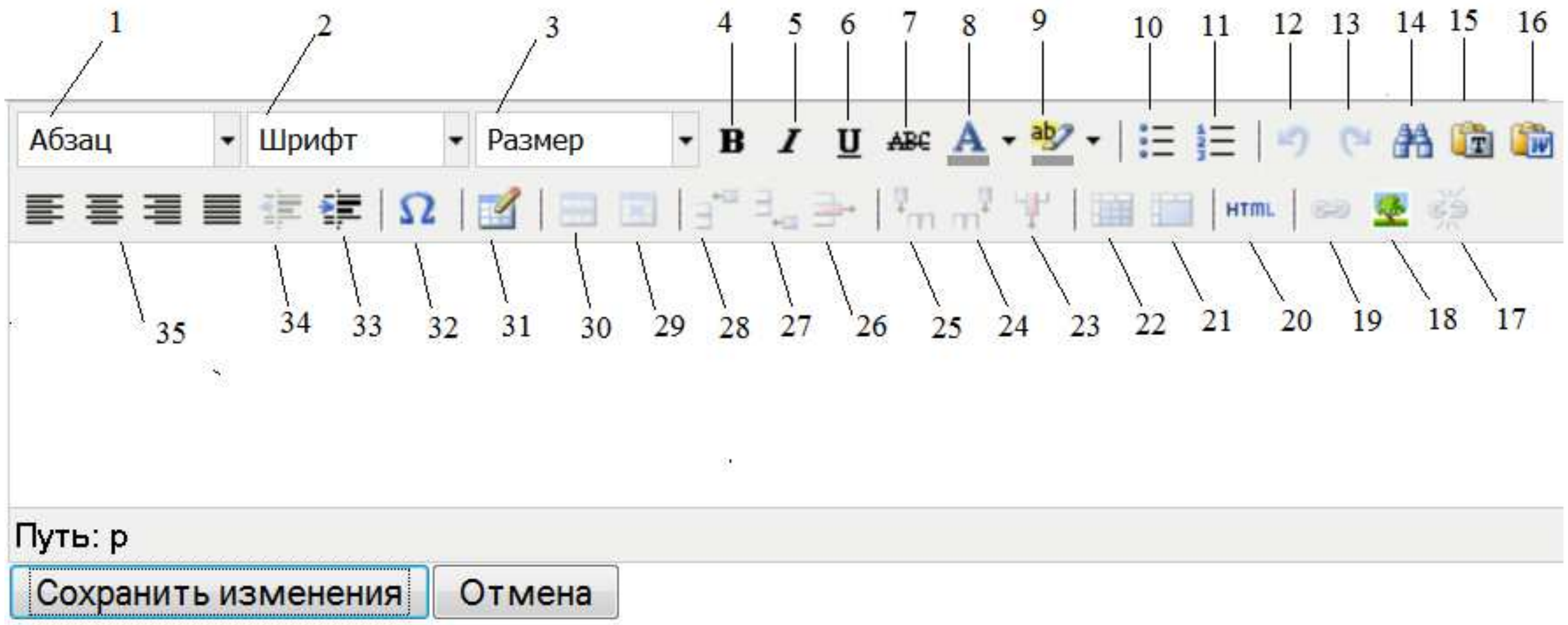

1 – список стилів;
2 – список типів шрифтів;
3 – список розмірів шрифтів;
4 – напівжирний;
5 – курсив;
6 – підкреслений;
7 – закреслений;
8 – вибрати колір тексту;
9 – вибрати колір фону;
10 – маркований список;
11 – нумерований список;
12 – відмінити;
13 – повернути;
14 – знайти;
15 – вставити неформатований текст;
16 – вставити з Word;
17 – видалити посилання;
18 – вставити/редагувати зображення;
19 – вставити/редагувати посилання;
20 – редагувати HTML код;
21 – об’єднати комірки;
22 – розділити комірки;
23 – видалити стовпчик;
24 – вставити стовпчик після;
25 – вставити стовпчик до;
26 – видалити рядок;
27 – вставити рядок після;
28 – вставити рядок до;
29 – властивості комірки;
30 – властивості рядка;
31 – вставити таблицю;
32 – вставити спеціальний символ;
33 – збільшити відступ;
34 – зменшити відступ;
35 – вирівнювання тексту.

*Рис. 8. Вікно редагування текстового редактору TinyMCE*

Зазначимо, що у вікні редагування, крім самого тексту, можна записувати вирази і формули в математичній нотації. Для цього необхідно скористатися командами мови LaTeX (рис. 9) і після натискання кнопки «Сохранить изменения» у полі виведення з’явиться вираз у природному математичному записі.

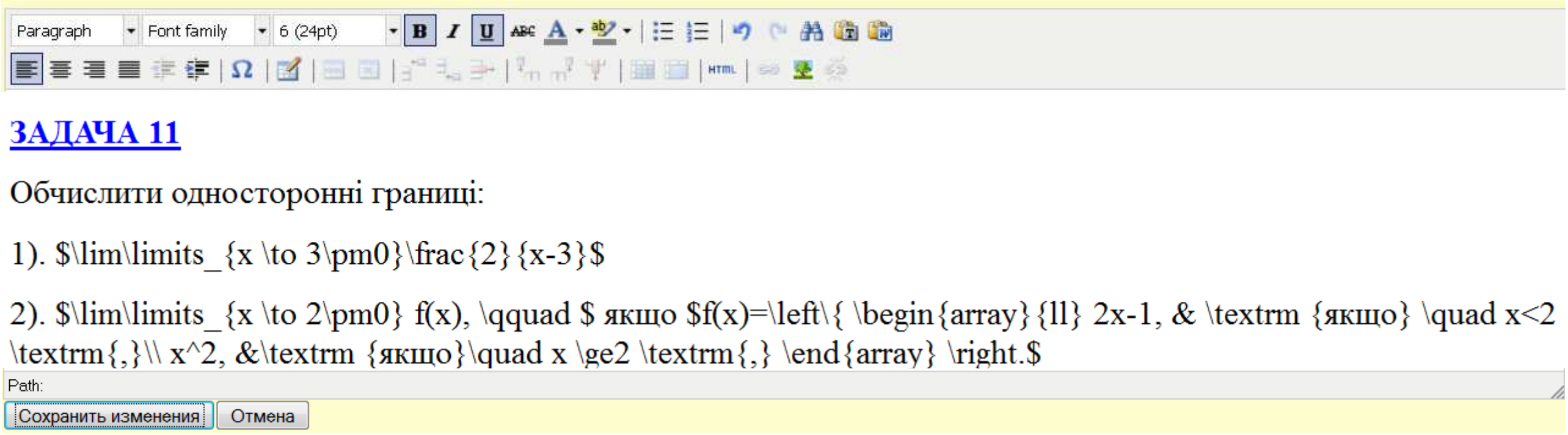

*а)*

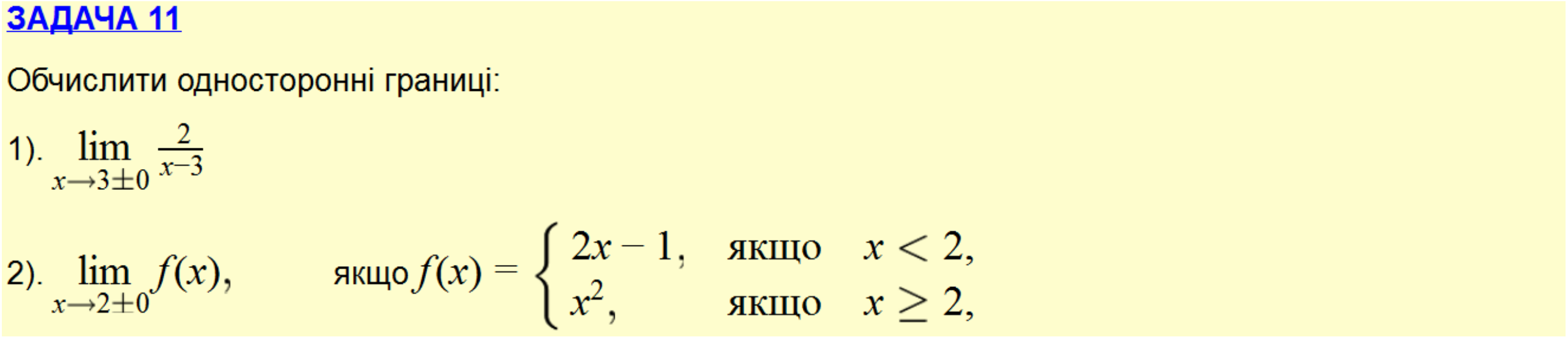

*б)*

*Рис. 9. Ілюстрація використання команд LaTeX для виведення математичного виразу*

Отже, робочі аркуші Sage поєднують в собі текстовий процесор, редактор формул і обчислювальні засоби, що надає можливість створювати у ММС інтерактивні математичні тексти.

**Висновки**. 1. Мобільне математичне середовище є інноваційним засобом навчання математичних дисциплін, визначальними особливостями якого є: об'єднання в собі інших засобів навчання (лекційних демонстрацій, динамічних моделей, тренажерів та навчальних експертних систем); можливість налаштування на конкретну математичну дисципліну; динамічна природа навчальних матеріалів.

2. Розробка навчально-методичної складової ММС «Вища математика» передбачає вміння: створювати елементи управління (за допомогою відповідних функцій обчислювального ядра − Web-СКМ Sage); описувати програмний код моделі (за допомогою мови програмування Python); додавати підписи різного формату у вигляді тексту або математичних виразів (за допомогою команд мови HTML і мови LaTeX); створювати інтерактивні математичні тексти (за допомогою текстового редактора TinyMCE); розробляти бази знань навчальних експертних систем (за допомогою eXpertise2Go).

3. Перспективними напрямами розвитку ММС є: розробка засобів доступу до ММС з мобільних пристроїв на платформах Google Android і Apple IOS;

вбудовування робочих аркушів ММС і їх компонентів у системи підтримки навчання; інтеграція із соціальними мережами; розробка засобів доступу до Wolfram|Alpha.

## Список використаних джерел

1. *Семеріков С. О.* Теорія та методика застосування мобільних математичних середовищ у процесі навчання вищої математики студентів економічних спеціальностей [Електронний ресурс] / Семеріков Сергій Олексійович, Словак Катерина Іванівна // Інформаційні технології і засоби навчання. – 2011. – №1(21). – Режим доступу до журналу : http://journal.iitta.gov.ua.

2. *Словак К. І.* Мобільне математичне середовище як новий засіб підвищення ефективності навчальної діяльності студентів з вищої математики / К. І. Словак // Інноваційні інформаційно-комунікаційні технології навчання математики, фізики, інформатики у середніх та вищих навчальних закладах : зб. наук. праць за матеріалами Всеукр. наук.-метод. конф. молодих науковців, 17–18 лютого 2011 р. – Кривий Ріг : Криворізький державний педагогічний ун-т, 2011. – С. 73–76.

3. Мобільне математичне середовище «Вища математика» [Електронний ресурс] / [К. І. Словак]. – 2011. – Режим доступу : http://korpus21.dyndns.org:8000/.

4. *Словак К.* І. Застосування мобільного математичного середовища SAGE у процесі навчання вищої математики студентів економічних ВНЗ / К. І. Словак // Педагогічні науки: теорія, історія, інноваційні технології : науковий журнал. – Суми : СумДПУ ім. А. С. Макаренка, 2010. – № 2 (4). – С. 345–354.

5. *Словак К. І.* Застосування ММС Sage у процесі навчання вищої математики / К. І. Словак // Вісник Черкаського університету. Серія педагогічні науки. – Вип. 191. – Ч. 1. – Черкаси : Вид. від. ЧНУ ім. Б. Хмельницького, 2010. – С. 106–111.

6. *Словак К. І.* Лекційні демонстрації у курсі вищої математики / К. І. Словак, М. В. Попель // Новітні комп’ютерні технології : матеріали VIII Міжнар. наук.-техн. конф. : К.–Севастополь, 14-17 вересня 2010 р. – К. : Міністерство регіонального розвитку та будівництва України, 2010. – С. 142–143.

7. *Словак К. І.* Мобільні математичні середовища: сучасний стан та перспективи розвитку / Словак К. І., Семеріков С. О., Триус Ю. В. // Науковий часопис Національного педагогічного університету імені М. П. Драгоманова. Серія №2. Комп’ютерно-орієнтовані системи навчання : зб. наук. праць / Редрада. – К. : НПУ імені М. П. Драгоманова, 2012. – № 12 (19). – С. 102–109.

8. *Михалевич В. М.* Математичні моделі генерування завдань з інтегрування частинами невизначених інтегралів / В. М. Михалевич, Я. В. Крупський, О. І. Шевчук // Вісник Вінницького політехнічного інституту. – 2008. – № 1. – С. 116–122.

9. Hawking L. e2gLite Tutorial [Electronic resource] / Louise Hawking – 2008. – 26 p. – Mode of access : http://ebookbrowse.com/e2glite-tutorial-pdf-d91750495.

10. Khan F. S. Dr. Wheat: A Web-based Expert System for Diagnosis of Dis-eases and Pest in Pakistani Wheat / Fohad Shahbaz Khan, Sead Razzaq, Kashif Irfan et al // Proceedings of The World Congress on Engineering 2008, July 2-4, 2008. – Vol I. – London, 2008. – P. 549–554.

11. Web-Enabled Expert System and Decision Table Software Demonstrations and Tutorials [Electronic resource] / eXpertise2Go.com. – 2009. – Mode of access : http://expertise2go.com

# МЕТОДИКА ПОСТРОЕНИЯ ОТДЕЛЬНЫХ КОМПОНЕНТОВ МОБИЛЬНОЙ МАТЕМАТИЧЕСКОЙ СРЕДЫ «ВЫСШАЯ МАТЕМАТИКА»

**Словак Екатерина Ивановна**, кандидат педагогических наук, старший преподаватель кафедры высшей математики, Криворожский экономический институт ГВУЗ «Криворожский национальный университет», г. Кривой Рог, e-mail: Slovak_kat@mail.ru

## Аннотация

Актуальность материала, изложенного в статье, обусловлена необходимостью разработки и внедрения в процесс обучения высокотехнологической информационно-коммуникационной образовательно-научной среды. В работе рассмотрен один из примеров такой web-ориентированной среды для обучения математическим дисциплинам студентов ВУЗов − мобильная математическая среда. На примере ММС «Высшая математика» продемонстрированы основные технологии и средства для построения таких компонентов среды, как лекционные демонстрации, динамические модели, тренажеры, генераторы учебных задач, учебные экспертные системы, индивидуальные домашние задания, примеры решения задач, задачи для практических занятий и т. п.

# METHODOLOGY OF SEPARATE COMPONENTS FORMATION OF MOBILE MATHEMATICAL ENVIRONMENT «HIGHER MATHEMATICS»

**Kateryna I. Slovak,** PhD (pedagogical sciences), senior lecturer of Department of higher mathematics, Kryvyi Rih Economic Institute of the State University « Kryvyi Rih National University», Kryvyi Rih, e-mail: Slovak_kat@mail.ru

**Resume**

The actuality of material stated in the article is caused by the necessity to develop and implement high-tech information and communication, educational and scientific environment to the leaning process. One of the examples of such web-based environment for learning mathematics of students of universities is discussed in the article. Basic technologies and means for forming such components of environment as lecture demonstrations, dynamic models, simulators, generators of leaning tasks, leaning expert systems, individual home tasks, examples of doing sums, tasks for practical lessons are demonstrated on the example of MME «HIGHER MATHEMATICS».